\documentclass[11pt]{amsart}

\usepackage{amsmath,amsthm, amscd, amssymb, amsfonts,mathtools}
\usepackage[all]{xy}
\usepackage{enumitem}

\usepackage{t1enc}
\usepackage[latin1]{inputenc}

\usepackage{fontsmpl}
\usepackage{srcltx}

\allowdisplaybreaks

\newcommand{\toba}{{\mathcal B}}

\newcommand{\ideal}{\mathcal M}

\newcommand{\vsp}{\vspace*{-0.7cm}}

\def\G{\mathbb{G}}

\newcommand{\ztu}{\overline{\zeta}}

\usepackage{color}

\numberwithin{equation}{section}\theoremstyle{plain}

\newtheorem{theorem}{Theorem}[section]

\newtheorem{prop}{Proposition}[subsection]

\theoremstyle{definition}

\theoremstyle{remark}
\newtheorem{rem}[theorem]{Remark}

\newcommand{\ydh}{{}^{H}_{H}\mathcal{YD}}

\def\G{\mathbb{G}}

\newcommand\ord{\operatorname{ord}}

\newcommand\Alg{\operatorname{Alg}}
\newcommand\Cleft{\operatorname{Cleft}}

\newcommand\gr{\operatorname{gr}}

\def\bq{\mathbf{q}}

\def\k{\Bbbk}

\def\ot{\otimes}

\def\Z{\mathbb{Z}}
\def\N{\mathbb{N}}

\def\B{\mathfrak{B}}

\def\eps{\epsilon}

\def\mH{\mathcal{H}}
\def\mL{\mathcal{L}}

\def\mJ{\mathcal{J}}

\def\mA{\mathcal{A}}
\def\mE{\mathcal{E}}
\def\mL{\mathcal{L}}

\def\br{\mathfrak{br}}

\def\ufo{\mathfrak{ufo}}

\def\hmod{H\text{-mod}}

\newcommand{\Gc}{{\mathcal G}}

\def\qb{\mathfrak{q}}

\def\lg{\langle}
\def\rg{\rangle}

\def\pf{\begin{proof}}
\def\epf{\end{proof}}

\newcommand{\Dchaintwo}[3]{\xymatrix@C-4pt{\overset{#1}{\underset{\ }{\circ}}\ar
@{-}[r]^{#2}
& \overset{#3}{\underset{\ }{\circ}}}}

\begin{document}


 \title[Modular and UFO liftings]{Examples of liftings of modular and unidentified type: $\ufo(7,8)$ and $\br(2,a)$}

\author[Garc\'ia Iglesias]{Agust\'in Garc\'ia Iglesias} \address{FaMAF-CIEM (CONICET), Universidad Nacional de C\'ordoba,
	Medina A\-llen\-de s/n, Ciudad Universitaria, 5000 C\' ordoba, Rep\'
	ublica Argentina.} \email{aigarcia@famaf.unc.edu.ar}

\author[Pacheco Rodr\'iguez]{Edwin Pacheco Rodr\'iguez} \address{Instituto de desarrollo económico e innovación (IDEI), Universidad Nacional de Tierra del Fuego (UNTDF), Thorne 302, 9420, Río Grande, Re\-pú\-bli\-ca Argentina} \email{efpacheco@untdf.edu.ar}

\dedicatory{To Nicol\'as Andruskiewitsch on his 60th birthday with gratitude.}

\thanks{\noindent 2000 \emph{Mathematics Subject Classification.}16W30. 
\newline The work was partially supported by CONICET, FONCyT-ANPCyT 2015-2845, Secyt (UNC), the MathAmSud project GR2HOPF}

\begin{abstract}
We compute all liftings of braidings of unidentified types $\ufo(7)$ and $\ufo(8)$. Some of these examples had been previously computed by Helbig.
We also present a list of examples of liftings of modular type $\br(2,a)$. 
\end{abstract}

\maketitle

\begin{center}{\textit{\tiny Watching a coast as it slips by the ship is like thinking about an enigma. There it is before you, smiling, frowning, inviting, grand, mean, insipid, or savage, and always mute with an air of whispering, ``Come and find out''.}}\end{center}
\vspace*{-.5cm}
\begin{flushright}
	{\sc \tiny  Conrad}
\end{flushright}

\section{Introduction}

This paper is a new contribution to the series initiated in \cite{AAG} with the computation of all liftings of finite-dimensional Nichols algebras of Cartan type $A$ and followed by \cite{GJ} for Cartan type $G_2$ and \cite{AnG-survey}, where a classification theorem was achieved and illustrated with the liftings of Cartan type $B_2$.

Let $A$ be a Hopf algebra such that its coradical $A_0=H$ is a Hopf subalgebra and the infinitesimal braiding $V\in\ydh$ of $A$ is a braided vector space of diagonal type with finite-dimensional Nichols algebra $\B(V)$. Then $A$ is a cocycle deformation of the graded Hopf algebra $\gr A=\B(V)\# H$ and an explicit, recursive, algorithm to recover $A$ from $\gr A$ is available. A particularly nice feature of this algorithm is that it works for any cosemisimple Hopf algebra $H$ fitting as the coradical of $A$. However, the deformations of $A$ can be extremely complicated from a combinatorial point of view. This was first observed in \cite{AAG} when the root of 1 was of order 3 and it was then made very clear for a subcase in \cite{GJ}, where a relation takes more than 10 pages.

It is unclear how to proceed at this point for any given Nichols algebra, as the computers we have at hand do not seem to be powerful enough to deal with every case.

\medbreak

In this paper, we focus  on liftings of type $\ufo(7)$ and $\ufo(8)$, which are still manageable. A remark is in place here. Helbig had already computed some of these liftings in \cite{Helbig}; with two minor restrictions that we are able to avoid here. On the one hand, the coradical in \cite{Helbig} is always assumed to be the group algebra of a finite abelian group, and it has already been said that our algorithm works for any $H$. On the other,  one family of liftings of type $\ufo(8)$ could not be fully described with the methods in loc.cit., and the algorithm shows that it is strong enough to give a full answer in this case. On top of this, our algorithm shows that the liftings obtained are cocycle deformations of $\B(V)\# H$ on each case.

We also compute some examples of modular type $\br(2,a)$, with parameters $\zeta\in\G_3'$ and $q\in \G_N'$, $N\neq 3$. The Nichols algebras in this case have different presentations according to whether $q=-1$, $q=-\zeta$, or $q$ is a generic root of 1. We compute all the liftings in the first two cases and give some families of examples on the third. The main obstruction on this last case is that we cannot produce an explicit formula for the coproduct of the $M$th power of the longest root, where $M=N$, $3N$ or $N/3$, according to the case. We provide examples for a fixed $N$ on all three cases; we conjecture that the shape of these liftings carries on to the generic case.

The core of the proofs is completely computational\footnote{See the log files at \texttt{https://www.famaf.unc.edu.ar/$\sim$aigarcia/publicaciones.htm}.}, we use the computer program \cite{GAP} together with the package \cite{GBNP} to follow the lines of the algorithm developed in \cite{AAG}; see also the survey \cite{AnG-survey}. However, the proofs are not automatic and require a deep insight on the combinatorial structure of the relations defining the Nichols algebras involved to reach the results.

\medbreak

Our main results on the description of the liftings on each case  can be summarized in Theorems \ref{thm:ufo} and \ref{thm:br} below. 
The isomorphism classes can be determined following  \cite[Theorem 3.9]{AAG}.

\begin{theorem}\label{thm:ufo}
	Let $H$ be a cosemisimple Hopf algebra and let $\mL$ be a Hopf algebra whose infinitesimal braiding $V\in\ydh$ is a principal Yetter-Drinfeld realization of a braided vector space of type $\ufo(7)$ or $\ufo(8)$.
	\begin{enumerate}
		\item If $V$ is of type $\ufo(7)$, then there are $\lambda_1,\lambda_2\in\k$ subject to \eqref{eqn:mu} such that $\mL\simeq \mL(\lambda_1,\lambda_2)$ as in Propositions \ref{pro:lifting-ufo(7)-a},  \ref{pro:lifting-ufo(7)-b} or \ref{pro:lifting-ufo(7)-c}, according to whether the braiding of $V$ is associated to the diagram \emph{(\ref{eq:dynkin-ufo(7)} a)}, \emph{(\ref{eq:dynkin-ufo(7)} b)} or \emph{(\ref{eq:dynkin-ufo(7)} c)}.
		\item If $V$ is of type $\ufo(8)$, then there are $\lambda_1,\lambda_2,\lambda_{12}\in\k$ subject to \eqref{eqn:mu} such that $\mL\simeq \mL(\lambda_1,\lambda_2,\lambda_{12})$ as in Propositions \ref{pro:lifting-ufo(8)-a} or   \ref{pro:lifting-ufo(8)-b} or there are $\lambda_1,\lambda_2\in\k$ subject to \eqref{eqn:mu} such that $\mL\simeq \mL(\lambda_1,\lambda_2)$ as in Proposition \ref{pro:lifting-ufo(8)-c}, according to whether the braiding of $V$ is associated to the diagram \emph{(\ref{eq:dynkin-ufo(8)} a)}, \emph{(\ref{eq:dynkin-ufo(8)} b)} or \emph{(\ref{eq:dynkin-ufo(8)} c)}.
	\end{enumerate}
\end{theorem}
\pf
See Sections \ref{sec:ufo7} and \ref{sec:ufo8}.
\epf

\begin{theorem}\label{thm:br}
	Let $H$ be a cosemisimple Hopf algebra and let $\mL$ be a Hopf algebra whose infinitesimal braiding $V\in\ydh$ is a principal Yetter-Drinfeld realization of a braided vector space of type $\br(2,a)$, with parameters $\zeta\in \G_3'$ and $q\notin \G_3'$.
	\begin{enumerate}
		\item If $q=-1$, respectively $q=-\zeta$, then there are $\lambda_1,\lambda_2,\lambda_3,\lambda_{112}\in\k$ subject to \eqref{eqn:mu} such that $\mL\simeq \mL(\lambda_1,\lambda_2,\lambda_3,\lambda_{112})$ as in Proposition \ref{pro:lifting-br(2,a)-q-1}, respectively Proposition \ref{pro:lifting-br(2,a)-a-q-zeta}.
		\item If $q\in \G_N'$, $N=4,6,12$, then there are $\lambda_1,\lambda_2,\lambda_{112}\in\k$ subject to \eqref{eqn:mu} such that $\mL\simeq \mL(\lambda_1,\lambda_2,\lambda_{112})$ as in Propositions \ref{pro:lifting-br(2,a)-a-q-generico-N4}, 
		\ref{pro:lifting-br(2,a)-a-q-generico-N6}, 
		\ref{pro:lifting-br(2,a)-a-q-generico-N12}, respectively.
	\end{enumerate}	
\end{theorem}
\pf
See Section \ref{sec:br}.
\epf

\section*{Acknowledgement}
We dedicate this article to Nicol\'as Andrukiewitsch, in occasion of his 60th birthday. We thank him for his generous guidance and for always encouraging us to face enigmas.

\section{Preliminaries}

We work over an algebraically closed field of characteristic zero. If $n\in\N$, we denote by $\G_n$ the set of roots of 1 of order $n$ and by $\G'_n\subset \G_n$ the subset of primitive roots. All algebras, unadorned tensor products, etc.~are assumed to be defined over $\k$. If $A$ is an algebra, we denote by $\Alg(A,\k)$ the set of algebra maps $A\to \k$. We denote by $G(A)$ the group of grouplike elements in a Hopf algebra $A$.

\medbreak 

Let $A$ be a Hopf algebra, with coradical $A_0=H$ being a Hopf subalgebra of $A$, and let $A_n=\wedge^n A_0$ denote the $n$th term of the coradical filtration of $A$, so $A_n=\{x\in A:\Delta(x)\in A_0\ot A+A\ot A_{n-1}\}$. In this setting, this is a Hopf algebra filtration and the associated graded algebra $\gr A=\oplus_{n\geq 0}A_n/A_{n-1}$, $A_{-1}=0$, is a graded Hopf algebra, which in turn splits as a smash product $\gr A=R\#H$, for certain graded Hopf algebra $R=\oplus R^n$ in the braided category $\ydh$ of Yetter-Drinfeld modules over $H$. This is a connected algebra, that is $R^0=\k$, and the subalgebra generated by the first component $V\coloneqq R^1$ coincides with the Nichols algebra $\B(V)$ of the braided vector space $V$. These observations are the cornerstone of the lifting method developed by Andruskiewitsch and Schneider which eventually led to a large and successful collection of classification results, to which this paper is our contribution. 

A Hopf algebra $A'$ is a cocycle deformation of $A$ if $A'=A$ as coalgebras and the multiplication $m'$ of $A'$ is obtained by a deformation $\sigma\ast m\ast\sigma^{-1}$ of the multiplication $m$ of $A$ by a Hopf 2-cocycle $\sigma\colon A\ot A\to \k$. We write $A_{\sigma}\coloneqq A'$.

Recall that a (right) cleft object for $A$ is a (right) comodule algebra $B$ with trivial coinvariants such that $B\simeq A$ as $A$-comodules. In particular, there is a convolution-invertible comodule isomorphism $\gamma\colon A\to B$, with $\gamma(1)=1$, called a ``section''. We shall denote by $\Cleft(A)$ the set of (isomorphism classes) of cleft objects for $A$. To such a pair $(A,B)$, there is an associated Hopf algebra $L=L(A,B)$ in such a way that $B$ becomes a $(L,H)$-bicleft object. Moreover, $L$ is a cocycle deformation of $A$, and every cocycle deformation $A_{\sigma}$ arises as $L(A_{(\sigma)},A)$ for some $A_{(\sigma)}\in\Cleft(A)$. The subscript ${}_{(\sigma)}$ not only indicates the connection to the Hopf algebra $A_{\sigma}$: the algebra structure on this cleft object arises itself as a (one-sided) deformation $\sigma\ast m$ of the multiplication of the algebra $A$. See \cite{S} for details on this construction and the associated equivalence.  

We refer the reader to \cite{Mo} for a full insight on Hopf algebra theory.

\subsection{The lifting method}

Let us briefly recall these ideas, which are also explained in detailed in \cite{AS-survey}. For a fixed cosemisimple Hopf algebra $H$ (originally $H=\k\Gamma$, $\Gamma$ a finite abelian group), the method proposes to compute all Yetter-Drinfeld modules $V\in \ydh$ with $\dim \B(V)<\infty$ and give a description of $\B(V)$ in terms of generators and relations, in order to lift or deform these relations to produce all Hopf algebras $A$ with $\gr A\simeq \B(V)\# H$. An independent and key step of this method, is to verify that any $A$ with $A_0\simeq H$ is ``generated in degree one'', that is that one always obtains $R=\B(V)$ in the decomposition $\gr A=R\# H$ mentioned above. Indeed, it was posed as a conjecture, still open in general, that this is aways the case.

When $H=\k\Gamma$ is finite and abelian, the braiding on any $V\in\ydh$ is diagonal; this is to say that there is basis $\{x_1,\dots,x_\theta\}$ of $V$ and a matrix $\bq=(q_{ij})\in \k^{\theta\times \theta}$ such that $c=c^{\qb}$ is determined by the formula
\[
c(x_i\ot x_j)=q_{ij}\,x_j\ot x_i, \qquad 1\leq i,j\leq \theta.
\]

This is the context of our work, and indeed the context in which the lifting method has been more fruitful, due to the work of many authors. All braided vector spaces with finite-dimensional Nichols algebras were classified by Heckenberger in \cite{H} in terms of generalized Dynkin diagrams (a finite collection of graphs with labels) codifying the matrices $\qb$. With this information, Angiono was able to give in \cite{A-jems-relations,Ang-crelle-gen_degree1} a complete presentation of the Nichols algebra attached to each such matrix, and used this presentation to confirm the conjecture about the generation in degree one in \cite{Ang-crelle-gen_degree1}. A big breakthrough was already present in \cite{AS-annals}, where a classification of the liftings was presented for the case $H=\k \Gamma$, with some restrictions on the order of the abelian group $\Gamma$. 

At this point, Masuoka showed in \cite{M} that every lifting defined in \cite{AS-annals} was a cocycle deformation of the associated graded Hopf algebra. This was the spark that initiated the program in \cite{AAGMV}, where an algorithm to produce these deformations as cocycle deformations of graded Hopf algebras was introduced, and later refined for the diagonal case in \cite{AAG}, and which led to the classification in \cite{AnG}. 

\medbreak 

Nonetheless, the explicit lifting of the relations is a work in progress that has to be approached on a case-by-case basis.

\medbreak

\subsection{Realizations}
We stress that the starting point of the lifting problem shifts from a fixed group $\Gamma$ to a fixed braided vector space $V$ with $\dim\B(V)<\infty$ such that $V\in\ydh$. In other words, one starts with a braided vector space (of diagonal type) -equivalently, with a braiding matrix $\qb=(q_{ij})_{1\leq i,j\leq\theta}$- and a cosemisimple Hopf algebra $H$ with a {\it principal Yetter-Drinfeld realization} $V\in\ydh$. This is encoded in the existence of grouplike elements $g_1,\dots,g_\theta\in G(H)$ and characters $\chi_1,\dots,\chi_\theta\in \Alg(H,\k)$ such that $\chi_j(g_i)=q_{ij}$, $1\leq i,j\leq\theta$.

\subsection{A hitchhiker's guide to the liftings}

A step-by-step description of the algorithm is presented in \cite{AnG-survey}, with the hope of recruiting volunteers to produce every lifting. We give a few highlights next.
The main ideas can be summarized in ``(a) deform one (primitive) relation at a time'' and ``(b) compute the cleft objects first''. The bedrock on which this philosophy stands is the following.
\begin{enumerate}[leftmargin=*]
	\item The ideal $\mJ(V)$ defining $\B(V)=T(V)/\mJ(V)$ is finitely generated by an homogeneous set $\Gc$ that admitis and stratification $\Gc=\Gc_0\sqcup \dots \sqcup\Gc_\ell$ in such a way that (the image of) every $r\in\Gc_{i}$ is primitive in $\B_i\coloneqq T(V)/\lg \Gc_0\sqcup \dots \sqcup\Gc_{i-1} \rg$.
	\item If $\mA_i$ is a (right) cleft object for $\mH_i\coloneqq \B_i\# H$, then there is a concrete procedure, inspired by \cite{Gunther}, to produce cleft objects $\mA_{i+1}$ for $\mH_{i+1}$.
	\item[(2')] Moreover, the (generic) definition of $\mA_i$ does not depend on $H$ and there is a cleft object $\mE_i$ for $\B_i$ (in $H$-mod) so that $\mA_i=\mE_i\# H$.
	\item Given $\mA_i$, it is possible to describe the left Hopf algebra $\mL_i\coloneqq L(\mA_i,\mH_i)$, see \cite{S}, by generators and relations.
\end{enumerate}
The outcome of this is a family of Hopf algebras $\mL_{\ell+1}=L(\mA_{\ell},\B(V)\#H)$ for every cleft object $\mA_{\ell+1}$ of $\mH_{\ell+1}=\B(V)\#H$. These are cocycle deformations of $\B(V)\#H$ and satisfy $\gr \mL_{\ell+1}\simeq \B(V)\# H$ by construction. One verifies that this is indeed an exhaustive collection of all the liftings.

We strongly suggest to go through the instructions in \cite{AnG-survey} for a detailed recipe to construct the liftings. We sketch here the main features.

\medbreak 
Assume we have reached step $i$, and fix $r\in \Gc_i$ (step 0 being $\mE_0\coloneqq T(V)$). We thus have a pre-Nichols algebra $\B_i$ with an associated Hopf algebra $\mH_i$, together with a (right) cleft object $\mA_i=\mE_i\# H$ and a section $\gamma_i\colon \mH_i\to \mA_i$. We also have the corresponding Hopf algebra $\mL_i=L(\mA_i,\mH_i)$, which coacts (on the left) on $\mA_i$; we call $\delta_i\colon \mA_i\to \mL_i\ot \mA_i$ this coaction.

Let $\sum a_j\alpha_j\in\Z^{\theta}$ be the degree of $r$; this defines
\[
g_r=\prod g_j^{a_{j}}\in G(H), \qquad \chi_r=\prod \chi_j^{a_j}\in\Alg(H,\k);
\]
a grouplike element and a character associated to $r$ and determined by the realization $V\in \ydh$.

\medbreak 
The algorithm runs as follows:
\begin{enumerate}
	\item Compute $r'\coloneqq \gamma_i(r)\in\mA_i$ and set 
	\[
	\mE_{i+1}=\mE_i/\lg r'-\lambda_r \rg,\]
	where 
	\begin{equation}\label{eqn:mu}
		\lambda_r=0 \qquad \text{if} \qquad \chi_r\neq \eps.
	\end{equation} 
	\item Compute $\tilde{r}\in\mL_i$ as the solution $\tilde{r}\ot 1_{\mA_i}=\delta_i(r)-g_{r}\ot r'$ and set
	\[
	\mL_{i+1}=\mL_i/\lg \tilde{r}-\lambda_r+\lambda_rg_{r} \rg.\]
	Notice that we may further normalize $\lambda_r$ as
	\begin{equation}\label{eqn:mu2}
		\lambda_r=0 \qquad \text{if} \qquad g_r=1.
	\end{equation} 
\end{enumerate}

\medbreak 

The computation of $r'$ and $\tilde{r}$ for each $r\in \Gc$ can be achieved in most cases with the help of the computer, but a prori it involves a deep understanding of the combinatorial structure of the Nichols algebra $\B(V)$.

\subsection*{Notation}

Fix $V\in\ydh$, with basis $\{x_1,\dots,x_\theta\}$; keep the notation as above. Notice both $\B(V)$ and $\mE$ are quotients of $T(V)$, while $\mL$ is a quotient of $T(V)\# H$ (as so are $\B(V)\# H$ and $\mA=\mE\# H$). In any case, we see that $\{x_1,\dots,x_\theta\}$ is a subset of generators. We shall rename these variables according to the context: we will keep the notation $\{x_1,\dots,x_\theta\}$ for the generators of $\B(V)$ (and $\mH$) and use $\{y_1,\dots,y_\theta\}$ for generators of $\mE$ (and $\mA$) and  $\{a_1,\dots,a_\theta\}$ for generators of $\mL$. As for the elements in $H$, particularly the group-likes $g_1,\dots,g_\theta$, we shall use this notation indistinctly.

\medbreak

We write the Dynkin diagrams corresponding to a given braiding in a single line, one next to the other, and refer to them as ``diagram a), b), ... '', counting from left to right. 

\medbreak 

We present the set $\Gc$ of relations defining each Nichols algebra in different lines, according to a chosen stratification $\Gc=\Gc_0\sqcup \dots \sqcup\Gc_\ell$, from top (relations in $\Gc_0$) to bottom (relations in $\Gc_\ell$).
In the cases under consideration in this article, we have that $\ell\in\{0,1,2\}$.

\section{Type $\ufo(7)$, $\zeta \in \G'_{12}$}\label{sec:ufo7}

The liftings of a braided vector space of this type have been  computed in 
\cite[Theorems 5.17 (1)--(3)]{Helbig}, for the case $H=\k\Gamma$, $\Gamma$ a a finite abelian group. In this section we generalize this result for liftings over any cosemisimple Hopf algebra $H$, and show that these are cocycle deformations of $\B(V)\#H$.

The Weyl groupoid has objects with generalized Dynkin diagrams of the following shapes:
\begin{align}\label{eq:dynkin-ufo(7)}
	\begin{aligned}
		&\Dchaintwo{-\overline{\zeta}^{\,2}}{-\zeta ^3}{-\zeta ^2}&
		&\Dchaintwo{-\overline{\zeta}^{\,2}}{\overline{\zeta}}{-1}&
		&\Dchaintwo{-\zeta ^3}{\zeta}{-1}&
		&\Dchaintwo{-\zeta ^2}{-\zeta }{-1}&
		&\Dchaintwo{-\zeta ^3}{-\overline{\zeta}}{-1}
	\end{aligned}
\end{align}

\begin{rem}
	Notice that diagram (\ref{eq:dynkin-ufo(7)} d) is of the shape of (\ref{eq:dynkin-ufo(7)} b) but with $\zeta^5$
	instead of $\zeta$. Idem for diagram (\ref{eq:dynkin-ufo(7)} e) with respect to (\ref{eq:dynkin-ufo(7)} c). Hence we are left with cases a), b) and c).
\end{rem}

\subsection{The generalized Dynkin diagram \emph{(\ref{eq:dynkin-ufo(7)}
		a)}}\label{subsec:ufo(7)-a}

\

The Nichols algebra $\toba(V)$ is generated by $x_1, x_2$ with defining relations
\begin{align*}
	x_1^3&=0; & x_2^3&=0; &  [x_1,x_{122}]_c -\frac{\zeta^{10}(1+\zeta)q_{12}}{1+\zeta^3}x_{12}^2&=0.
\end{align*}

\begin{prop}\label{pro:cleft-ufo(7)-a}
	Fix $V$ of type $\ufo(7)$,
	as in \emph{(\ref{eq:dynkin-ufo(7)} a)} and $\lambda_1, \lambda_2\in\k$ subject to \eqref{eqn:mu}. 
	Let $\mE=\mE(\lambda_1,\lambda_2)$ be the quotient of $T(V)$ modulo the ideal generated by
	\begin{align*}
		y_1^3&=\lambda_1; & y_2^3&=\lambda_2;  &  [y_1,& y_{122}]_c
		-\frac{\zeta^{10}(1+\zeta)q_{12}}{1+\zeta^3}y_{12}^2=0.
	\end{align*}
	Then $\mE\in\hmod$ and $\mA(\lambda_1,\lambda_2)\coloneqq \mE\#H \in\Cleft\B(V)\#H$.
\end{prop}
\pf
The relations defining $\toba(V)$ are all primitive. Notice that the leftmost relation has degree $2\alpha_1+2\alpha_2$ and cannot be deformed by \eqref{eqn:mu}. Indeed, if $\chi_1^2\chi_2^2=\eps$, then $q_{11}^2q_{12}=\zeta^{-4}q_{12}^2=1$ and $q_{21}^2q_{22}=q_{12}^2\zeta^{4}=1$, which gives $\zeta^6=(q_{12}q_{21})^2=1$, a contradiction since $\zeta\in \G'_{12}$.

On the other hand, $\lambda_1\neq 0$ only if $q_{21}^3=1$ and $\lambda_2\neq 0$ only if $q_{12}^3=1$. In particular, we cannot have $\lambda_1\lambda_2\neq 0$, as otherwise $(-\zeta^3)^3=(q_{12}q_{21})^3=1$, again a contradiction with $\zeta\in \G'_{12}$.
\epf

\begin{prop}\label{pro:lifting-ufo(7)-a}
	Fix $V$ of type $\ufo(7)$, as in \emph{(\ref{eq:dynkin-ufo(7)} a)},  and $\lambda_1, \lambda_2\in\k$ subject to \eqref{eqn:mu}. Let $\mL=\mL(\lambda_1,\lambda_2)$ be the quotient 
	$T(V)\# H / \ideal$, where
	$\ideal$ is the ideal of $T(V)\# H$ generated by
	\begin{align*}
		a_1^3&=\lambda_1(1-g_1^3); & a_2^3&=\lambda_2(1-g_2^3); \\
		[a_1,& a_{122}]_c
		-\frac{\zeta^{10}(1+\zeta)q_{12}}{1+\zeta^3}a_{12}^2=0.
	\end{align*}
	Then $\mL$ is a Hopf algebra and a cocycle deformation of $\B(V)\# H$ such that $\gr\mL\simeq \B(V)\# H$.
	
	Conversely, if $L$ is a Hopf algebra such that its infinitesimal braiding is a principal realization $V\in\ydh$, then there are $\lambda_1, \lambda_2\in\k$ such that 
	$L\simeq \mL(\lambda_1, \lambda_2)$.
\end{prop}
\pf
The relations defining $\toba(V)$ are all primitive.
\epf

\subsection{The generalized Dynkin diagram \emph{(\ref{eq:dynkin-ufo(7)}
		b)}}\label{subsec:ufo(7)-b}

\

The Nichols algebra $\toba(V)$ is generated by $x_1, x_2$ with defining relations
\begin{align*}
	x_1^3&=0; & x_2^2&=0; 
\end{align*}
\vsp
\begin{align*}
	[[x_{112},x_{12}]_c,x_{12}]_c&=0.
\end{align*}

\begin{prop}\label{pro:cleft-ufo(7)-b}
	Fix $V$ of type
	$\ufo(7)$,
	as in \emph{(\ref{eq:dynkin-ufo(7)} b)} and $\lambda_1, \lambda_2\in\k$ subject to \eqref{eqn:mu}. Let $\mE=\mE(\lambda_1,\lambda_2)$ be the quotient of $T(V)$ modulo the ideal generated by
	\begin{align*}
		y_1^3&=\lambda_1; & y_2^2&=\lambda_2; & [[y_{112}&,y_{12}]_c,y_{12}]_c=0.
	\end{align*}
	Then $\mE\in\hmod$ and $\mA(\lambda_1,\lambda_2)\coloneqq \mE\#H \in\Cleft\B(V)\#H$.
\end{prop}
\pf
See \texttt{cleft-ufo7b.(g|log)}. As in (the proof of) Proposition \ref{pro:cleft-ufo(7)-a}, we have $\lambda_1\lambda_2=0$.
\epf

\begin{prop}\label{pro:lifting-ufo(7)-b}
	Fix $V$ of type $\ufo(7)$, as in \emph{(\ref{eq:dynkin-ufo(7)} b)},  and $\lambda_1, \lambda_2\in\k$ subject to \eqref{eqn:mu}. Let $\mL=\mL(\lambda_1,\lambda_2)$ be the quotient 
	$T(V)\# H / \ideal$, where
	$\ideal$ is the ideal of $T(V)\# H$ generated by
	\begin{align*}
		a_1^3&=\lambda_1(1-g_1^3); \qquad a_2^2=\lambda_2(1-g_2^2); \\ 
		[[a_{112}&,a_{12}]_c,a_{12}]_c=\lambda_2q_{12}(1+\zeta^7)a_{112}a_1^2g_2^2.
	\end{align*}
	Then $\mL$ is a Hopf algebra and a cocycle deformation of $\B(V)\# H$ such that $\gr\mL\simeq \B(V)\# H$.
	
	Conversely, if $L$ is a Hopf algebra such that its infinitesimal braiding is a principal realization $V\in\ydh$, then there are $\lambda_1, \lambda_2\in\k$ such that 
	$L\simeq \mL(\lambda_1, \lambda_2)$.
\end{prop}
\pf
See \texttt{lift-ufo7b.(g|log)}.
\epf

\subsection{The generalized Dynkin diagram \emph{(\ref{eq:dynkin-ufo(7)}
		c)}}\label{subsec:ufo(7)-c}

\

The Nichols algebra $\toba(V)$ is generated by $x_1, x_2$ with defining relations
\begin{align*}
	x_1^4&=0; & x_2^2&=0;
\end{align*}
\vsp
\begin{align*}
	[ x_{112},x_{12}]_c&=0.
\end{align*}

\begin{prop}\label{pro:cleft-ufo(7)-c}
	Fix $V$ of type
	$\ufo(7)$,
	as in \emph{(\ref{eq:dynkin-ufo(7)} c)} and $\lambda_1, \lambda_2\in\k$ subject to \eqref{eqn:mu}. Let $\mE=\mE(\lambda_1,\lambda_2)$ be the quotient of $T(V)$ modulo the ideal generated by
	\begin{align*}
		y_1^4&=\lambda_1 & y_2^2&=\lambda_2, &  [ y_{112},y_{12}]_c&=0.
	\end{align*}
	Then $\mE\in\hmod$ and $\mA(\lambda_1,\lambda_2)\coloneqq \mE\#H \in\Cleft\B(V)\#H$.
\end{prop}
\pf
See \texttt{cleft-ufo7c.(g|log)}. Again, it follows that $\lambda_1\lambda_2=0$.
\epf

\begin{prop}\label{pro:lifting-ufo(7)-c}
	Fix $V$ of type $\ufo(7)$, as in \emph{(\ref{eq:dynkin-ufo(7)} c)},  and $\lambda_1, \lambda_2\in\k$ subject to \eqref{eqn:mu}. Let $\mL=\mL(\lambda_1,\lambda_2)$ be the quotient 
	$T(V)\# H / \ideal$, where
	$\ideal$ is the ideal of $T(V)\# H$ generated by
	\begin{align*}
		a_1^4=\lambda_1(1-g_1^4); \qquad a_2^2=\lambda_2(1-g_2^2); \\
		[ a_{112},a_{12}]_c=\lambda_2q_{12}(-2\zeta^7+2\zeta^8-\zeta^{11})g_2^2a_1^3.
	\end{align*}
	Then $\mL$ is a Hopf algebra and a cocycle deformation of $\B(V)\# H$ such that $\gr\mL\simeq \B(V)\# H$.
	
	Conversely, if $L$ is a Hopf algebra such that its infinitesimal braiding is a principal realization $V\in\ydh$, then there are $\lambda_1, \lambda_2\in\k$ such that 
	$L\simeq \mL(\lambda_1, \lambda_2)$.
\end{prop}

\pf
See \texttt{lift-ufo7c.(g|log)}.
\epf

\section{Type $\ufo(8)$, $\zeta \in \G'_{12}$}\label{sec:ufo8}
The Weyl groupoid has objects with generalized Dynkin diagrams of the following shapes:
\begin{align}\label{eq:dynkin-ufo(8)}
	\begin{aligned}
		&\Dchaintwo{-\zeta ^2}{\zeta }{-\zeta ^2}
		& &\Dchaintwo{-\zeta ^2}{\zeta ^3}{-1}
		& & \Dchaintwo{-\ztu}{-\zeta ^3}{-1}
	\end{aligned}
\end{align}

\subsection{The generalized Dynkin diagram \emph{(\ref{eq:dynkin-ufo(8)}
		a)}}\label{subsec:ufo(8)-a}

\

The Nichols algebra $\toba(V)$ is generated by $x_1, x_2$ with defining relations
\begin{align*}
	x_1^3&=0; & x_2^3&=0; & [x_1,x_{122}]_c-(1+\zeta+\zeta^2)\zeta^4q_{12}x_{12}^2=0;
\end{align*}
\vsp
\begin{align*}
	x_{12}^{12}&=0.
\end{align*}

\begin{prop}\label{pro:cleft-ufo(8)-a}
	Fix $V$ of type
	$\ufo(8)$,
	as in \emph{(\ref{eq:dynkin-ufo(8)} a)} and $\lambda_1,\lambda_2,\lambda_{12}\in\k$ subject to \eqref{eqn:mu}. Let $\mE=\mE(\lambda_1,\lambda_2,\lambda_{12})$ be the quotient of $T(V)$ modulo the ideal generated by
	\begin{align*}
		y_1^3&=\lambda_1; & y_2^3&=\lambda_2; &
		y_{12}^{12}&=\lambda_{12};
	\end{align*}
	\vsp
	\begin{align*}
		[y_1,y_{122}]_c-(1+\zeta+\zeta^2)\zeta^4q_{12}y_{12}^2=0;
	\end{align*}
	Then $\mE\in\hmod$ and $\mA(\lambda_1,\lambda_2,\lambda_{12})\coloneqq \mE\#H \in\Cleft\B(V)\#H$.
\end{prop}

\pf
See \texttt{cleft-ufo8a.(g|log)}.

Notice that we cannot have $\chi_1^2\chi_2^2=\eps$, hence the relation of degree $2\alpha_1+2\alpha_2$ is not deformed, by \eqref{eqn:mu}. On the other hand, $\lambda_1\neq 0$ only if $q_{21}^3=1$ and $\lambda_2\neq 0$ only if $q_{12}^3=1$. In particular, $\lambda_1\lambda_2=0$.
\epf

\begin{prop}\label{pro:lifting-ufo(8)-a}
	Fix $V$ of type
	$\ufo(8)$,
	as in \emph{(\ref{eq:dynkin-ufo(8)} a)},  and $\lambda_1, \lambda_2,\lambda_{12}\in\k$ subject to \eqref{eqn:mu}. Let $\mL=\mL(\lambda_1,\lambda_2,\lambda_{12})$ be the quotient 
	$T(V)\# H / \ideal$, where
	$\ideal$ is the ideal of $T(V)\# H$ generated by
	\begin{align*}
		a_1^3&=\lambda_1(1-g_1^3); & a_2^3&=\lambda_2(1-g_2^3); &
		a_{12}^{12}&=\lambda_{12}(1-g_1^{12}g_2^{12});
	\end{align*}
	\vsp
	\begin{align*}
		[a_1,a_{122}]_c-(1+\zeta+\zeta^2)\zeta^4q_{12}a_{12}^2=0;
	\end{align*}
	Then $\mL$ is a Hopf algebra and a cocycle deformation of $\B(V)\# H$ such that $\gr\mL\simeq \B(V)\# H$.
	
	Conversely, if $L$ is a Hopf algebra such that its infinitesimal braiding is a principal realization $V\in\ydh$, then there are $\lambda_1, \lambda_2,\lambda_{12}\in\k$ such that 
	$L\simeq \mL(\lambda_1, \lambda_2,\lambda_{12})$.
\end{prop}
\pf
See \texttt{lift-ufo8a.(g|log)}.
\epf

\subsection{The generalized Dynkin diagram \emph{(\ref{eq:dynkin-ufo(8)}
		b)}}\label{subsec:ufo(8)-b}

\

The Nichols algebra $\toba(V)$ is generated by $x_1, x_2$ with defining relations
\begin{align*}
	x_1^3&=0; & x_2^2&=0; 
\end{align*}
\vsp
\begin{align*}
	[[x_{112},x_{12}]_c,x_{12}]_c&=0;
\end{align*}
\vsp
\begin{align*}
	x_{12}^{12}&=0.
\end{align*}

\begin{prop}\label{pro:cleft-ufo(8)-b}
	Fix $V$ of type
	$\ufo(8)$,
	as in \emph{(\ref{eq:dynkin-ufo(8)} b)} and $\lambda_1, \lambda_2,\lambda_{12}\in\k$ subject to \eqref{eqn:mu}.  Let $\mE=\mE(\lambda_1,\lambda_2,\lambda_{12})$ be the quotient of $T(V)$ modulo the ideal generated by
	\begin{align*}
		y_1^3&=\lambda_1; & y_2^2&=\lambda_2; & [[y_{112},& y_{12}]_c,y_{12}]_c=0 & 
		y_{12}^{12}&=\lambda_{12}; 
	\end{align*}
	Then $\mE\in\hmod$ and $\mA(\lambda_1,\lambda_2,\lambda_{12})\coloneqq \mE\#H \in\Cleft\B(V)\#H$.
\end{prop}

\pf
See \texttt{cleft-ufo8b-1.(g|log)} and \texttt{cleft-ufo8b-2.(g|log)}. Here, as above, the relation of degree $4\alpha_1+3\alpha_2$ cannot be deformed and $\lambda_1\lambda_2=0$, by \eqref{eqn:mu}.
\epf

\begin{prop}\label{pro:lifting-ufo(8)-b}
	Fix $V$ of type
	$\ufo(8)$,
	as in \emph{(\ref{eq:dynkin-ufo(8)} b)},  and $\lambda_1, \lambda_2,\lambda_{12}\in\k$ subject to \eqref{eqn:mu}. Let $\mL=\mL(\lambda_1,\lambda_2,\lambda_{12})$ be the quotient 
	$T(V)\# H / \ideal$, where
	$\ideal$ is the ideal of $T(V)\# H$ generated by
	\begin{align*}
		a_1^3&=\lambda_1(1-g_1^3); \qquad  a_2^2=\lambda_2(1-g_2^2); \\
		[[a_{112}, a_{12}]_c,a_{12}]_c&=-2q_{12}\lambda_2(1+\zeta^3)g_2^2a_{112}a_1^2\\
		a_{12}^{12}&=\lambda_{12}(1-g_1^{12}g_2^{12}); 
	\end{align*}
	Then $\mL$ is a Hopf algebra and a cocycle deformation of $\B(V)\# H$ such that $\gr\mL\simeq \B(V)\# H$.
	
	Conversely, if $L$ is a Hopf algebra such that its infinitesimal braiding is a principal realization $V\in\ydh$, then there are $\lambda_1, \lambda_2,\lambda_{12}\in\k$ such that 
	$L\simeq \mL(\lambda_1, \lambda_2,\lambda_{12})$.
\end{prop}
\pf
See \texttt{lift-ufo8b-1.(g|log)} and \texttt{lift-ufo8b-2.(g|log)}. We remark that the relation of degree $4\alpha_1+3\alpha_2$, which does not produce a deformation parameter by itself by \eqref{eqn:mu}, is nevertheless deformed as a consequence of the deformation of relation $x_2^2=0$, i.e.~when $\lambda_2\neq 0$.
\epf

The liftings in Proposition \ref{pro:lifting-ufo(8)-b} have been previously 
investigated in \cite[Theorem 5.17 (4)]{Helbig},  for the case $H=\k\Gamma$, $\Gamma$ an abelian group. 
In this case, however, a complete description was not attained. Indeed, in \cite[Theorem 5.17 (4)(b)]{Helbig}, the deformation of the relation $x_{12}^{12}$ corresponding to diagram \emph{(\ref{eq:dynkin-ufo(8)} b)} is presented as $a_{12}^{12}=d_{12}$, for a certain unspecified element $d_{12}$. In particular, we show that $d_{12}=\lambda_{12}(1-g_1^{12}g_2^{12})$.

\subsection{The generalized Dynkin diagram \emph{(\ref{eq:dynkin-ufo(8)}
		c)}}\label{subsec:ufo(8)-c}

\

The Nichols algebra $\toba(V)$ is generated by $x_1, x_2$ with defining relations
\begin{align*}
	x_1^{12}&=0; & x_2^2&=0; & x_{11112}&=0; 
\end{align*}
\vsp
\begin{align*}
	[x_{112},x_{12}]_c&=0.
\end{align*}

\begin{prop}\label{pro:cleft-ufo(8)-c}
	Fix $V$ of type
	$\ufo(8)$,
	as in \emph{(\ref{eq:dynkin-ufo(8)} c)} and $\lambda_1, \lambda_2\in\k$ subject to \eqref{eqn:mu}.  Let $\mE=\mE(\lambda_1,\lambda_2)$ be the quotient of $T(V)$ modulo the ideal generated by
	\begin{align*}
		y_1^{12}&=\lambda_1; & y_2^2&=\lambda_2; &   y_{11112}&=0; &  [y_{112},y_{12}]_c&=0.
	\end{align*}
	Then $\mE\in\hmod$ and $\mA(\lambda_1,\lambda_2)\coloneqq \mE\#H \in\Cleft\B(V)\#H$.
\end{prop}

\pf
See \texttt{cleft-ufo8c.(g|log)}. The relations of degree $4\alpha_1+\alpha_2$ (the quantum Serre relation) and $3\alpha_1+2\alpha_2$ cannot be deformed. As well, $\lambda_1\neq 0$ only if $q_{21}^{12}=1$ and $\lambda_2\neq 0$ only if $q_{12}^2=1$. Hence, we can have $\lambda_1\lambda_2\neq 0$ if (and only if) $q_{12}=\pm1$ and $q_{21}=\mp\zeta^3$.
\epf

\begin{prop}\label{pro:lifting-ufo(8)-c}
	Fix $V$ of type $\ufo(8)$, as in \emph{(\ref{eq:dynkin-ufo(8)} c)},  and $\lambda_1, \lambda_2\in\k$ subject to \eqref{eqn:mu}. Let $\mL=\mL(\lambda_1,\lambda_2)$ be the quotient 
	$T(V)\# H / \ideal$, where
	$\ideal$ is the ideal of $T(V)\# H$ generated by
	\begin{align*}
		a_1^{12}&=\lambda_1(1-g_1^{12}); & a_2^2&=\lambda_2(1-g_2^2); \\
		a_{11112}&=0; &  [a_{112},a_{12}]_c&=2q_{12}\lambda_2g_2^2a_1^3;
	\end{align*}
	Then $\mL$ is a Hopf algebra and a cocycle deformation of $\B(V)\# H$ such that $\gr\mL\simeq \B(V)\# H$.
	
	Conversely, if $L$ is a Hopf algebra such that its infinitesimal braiding is a principal realization $V\in\ydh$, then there are $\lambda_1, \lambda_2\in\k$ such that 
	$L\simeq \mL(\lambda_1, \lambda_2)$.
\end{prop}
\pf
See \texttt{lift-ufo8c.(g|log)}. As above, we remark that the relation of degree $3\alpha_1+2\alpha_2$ does not produce a deformation parameter by itself by \eqref{eqn:mu} but is deformed when $\lambda_2\neq 0$.
\epf

The liftings in Proposition \ref{pro:lifting-ufo(8)-c} have been previously computed in 
\cite[Theorem 5.17 (5)]{Helbig},  for the case $H=\k\Gamma$, $\Gamma$ an abelian group.

\section{Type $\br(2, a)$, $\zeta \in \G_3$, $q\notin \G_3$}\label{sec:br}
The Weyl groupoid has objects with generalized Dynkin diagrams of the following shapes:
\begin{align}\label{eq:dynkin-br(2,a)}
	&\xymatrix{ \overset{\zeta}{\underset{\ }{\circ}} \ar  @{-}[r]^{q^{-1}}  &
		\overset{q}{\underset{\ }{\circ}}}
	& &\xymatrix{ \overset{\zeta}{\underset{\ }{\circ}} \ar  @{-}[r]^{\zeta^2q}  &
		\overset{\zeta q^{-1}}{\underset{\ }{\circ}}}
\end{align}
It is standard of type $B_2$ with Cartan matrix $\begin{pmatrix}2 & -2 \\-1 & 2
\end{pmatrix}$.

The second, rightmost, diagram is equivalent to the first one, with $\zeta q^{-1}$ instead of $q$. We shall restrict ourselves to the leftmost one.

We shall use the notation $g_{112}\coloneqq g_1^2g_2$, $g_{122}\coloneqq g_1g_2^2$; also $g_{11212}\coloneqq g_1^3g_2^2$.

\subsection{The case $q= -1$}\label{subsubsec:br(2,a)-q-1}
The Nichols
algebra $\toba(V)$ is generated by $x_1, x_{2}$
with defining relations
\begin{align*}
	x_{1}^3 &= 0; & x_{2}^{2} &=0; 
\end{align*}
\vsp
\begin{align*}
	[x_{112}, x_{12}]_c&=0; 
\end{align*}
\vsp
\begin{align*}
	x_{112}^{6} &=0.
\end{align*}

\begin{prop}\label{pro:cleft-br(2,a)-q-1}
	Fix $V$ of type
	$\br(2, a)$,
	as in \emph{(\ref{eq:dynkin-br(2,a)} a)}, with $q=-1$.  Let $\lambda_1,\lambda_2,\lambda_3,\lambda_{112}\in\k$ subject to \eqref{eqn:mu} and let $\mE=\mE(\lambda_1,\lambda_2,\lambda_3,\lambda_{112})$ be the quotient of $T(V)$ modulo the ideal generated by
	\begin{align*}
		y_{1}^3 &= \lambda_1; &  y_{2}^{2} &=\lambda_2; & [y_{112}, y_{12}]_c&=\lambda_3; & y_{112}^{6} &=\lambda_{112}.
	\end{align*}
	Then $\mE\in\hmod$ and $\mA(\lambda_1,\lambda_2,\lambda_3,\lambda_{112})\coloneqq \mE\#H \in\Cleft\B(V)\#H$.
\end{prop}

\pf
See \texttt{cleft-br2a-q-1.(g|log)}. In this case $\lambda_1\neq 0$ only if $q_{21}^3=1$ and $\lambda_2\neq 0$ only if $q_{12}^2=1$. For $\lambda_3\neq 0$, we need 
$\bq=\left(\begin{smallmatrix}\zeta& -1 \\ 1& -1
\end{smallmatrix}\right)$, which allows (and it is equivalent to) $\lambda_1\lambda_2\neq0$. 
\epf

\begin{prop}\label{pro:lifting-br(2,a)-q-1}
	Fix $V$ of type
	$\br(2, a)$,
	as in \emph{(\ref{eq:dynkin-br(2,a)} a)}, with $q=-1$. Let $\lambda_1, \lambda_2,\lambda_3,\lambda_{112}\in\k$ subject to \eqref{eqn:mu} and let $\mL=\mL(\lambda_1, \lambda_2,\lambda_3,\lambda_{112})$ be the quotient $T(V)\# H / \ideal$, where
	$\ideal$ is the ideal of $T(V)\# H$ generated by
	\begin{align*}
		a_{1}^3 &= \lambda_1(1-g_1^3); \qquad a_{2}^{2} =\lambda_2(1-g_2^2); \\
		[a_{112}, a_{12}]_c&=\lambda_3(1-g_1^3g_2^2)-4\lambda_1\lambda_2g_2^2(1-g_1^3);  \\
		a_{112}^{6} &=\lambda_{112}(1-g_{112}^6)-s_{112}(a)-s_{112}(g),
	\end{align*}
	where $g_{112}=g_1^2g_2$ and for $\omega=\zeta^2-\zeta$ and $\lambda_{123}=4\lambda_1\lambda_2-\lambda_3$:
	\begin{align*}
		s_{112}(a)&=2\lambda_1\lambda_{123}\omega\, g_{11212}a_{12}^2a_{112}^2 - 2\lambda_1\lambda_{123}\omega\,g_{11212}a_{2}a_{112}^3 \\
		&+ 8\lambda_1^2\lambda_2\lambda_{123}\omega\, g_{2}^2g_{11212}a_{12}a_{112} - 2\lambda_1\lambda_{123}^2\omega\, g_{11212}^2a_{12}a_{112} \\
		&- 2\lambda_1\lambda_3\lambda_{123}\omega\, g_{11212}a_{12}a_{112}; \\
		s_{112}(g)&=8\lambda_1^2\lambda_2(8\lambda_1^2\lambda_2^2-8\lambda_1\lambda_2\lambda_3+3\lambda_3^2)g_{2}^2 + 64\lambda_1^3\lambda_2^2(\lambda_1\lambda_2-\lambda_3)\, g_{2}^4\\
		&+ 64\lambda_1^4\lambda_2^3\, g_{2}^6 + \lambda_1\lambda_3^2\lambda_{123}\omega\, g_{11212} - 8\lambda_1^2\lambda_2\lambda_3\lambda_{123}\omega\, g_{2}^2g_{11212} \\
		&+ 16\lambda_1^3\lambda_2^2\lambda_{123}\omega\,g_{2}^4g_{11212} - 20\lambda_1^2\lambda_2\lambda_{123}^2g_{2}^2g_{11212}^2 \\
		&+ 4\lambda_1(16\lambda_1^3\lambda_2^3+8\lambda_1^2\lambda_2^2\lambda_3-7\lambda_1\lambda_2\lambda_3^2+\lambda_3^3)g_{11212}^2\\
		& +4\lambda_1(16\lambda_1^3\lambda_2^3-16\lambda_1^2\lambda_2^2\lambda_3+6\lambda_1\lambda_2\lambda_3^2-\lambda_3^3)g_{1}^3g_{11212}^3 \\
		&+ \lambda_1\omega(-64\lambda_1^3\lambda_2^3+48\lambda_1^2\lambda_2^2\lambda_3-12\lambda_1\lambda_2\lambda_3^2+\lambda_3^3)g_{11212}^3.
	\end{align*}

	Then $\mL$ is a Hopf algebra and a cocycle deformation of $\B(V)\# H$ such that $\gr\mL\simeq \B(V)\# H$.
	
	Conversely, if $L$ is a Hopf algebra such that its infinitesimal braiding is a principal realization $V\in\ydh$, then there are $\lambda_1, \lambda_2,\lambda_3,\lambda_{112}\in\k$ such that 
	$L\simeq \mL(\lambda_1, \lambda_2,\lambda_3,\lambda_{112})$.
	
\end{prop}

\pf 
See \texttt{lift-br2a-q-1-PBW-wg11212.(g|log)}. In this case, we use variables corresponding to the (induced) PBW basis of the liftings, and, interestingly, we need to add a variable for the element $g_{11212}=g_1^3g_2^2$ to get an expression of $s_{112}$ in PBW form.
\epf

\subsection{The case $q=-\zeta$}\label{subsubsec:br(2,a)-q-zeta}
$\toba(V)$ is
generated by $x_1, x_{2}$
with defining relations
\begin{align*}	
	x_{1}^3 &= 0; & x_{2}^{6} &=0; & x_{221} &=0;
\end{align*}
\vsp
\begin{align*}
	x_{112}^{2} &=0.
\end{align*}

\begin{prop}\label{pro:cleft-br(2,a)-q-zeta}
	Fix $V$ of type
	$\br(2, a)$,
	as in \emph{(\ref{eq:dynkin-br(2,a)} a)}, with $q=-\zeta$.
	Let $\lambda_1,\lambda_2,\lambda_{3},\lambda_{112}\in\k$ subject to \eqref{eqn:mu} and let $\mE=\mE(\lambda_1,\lambda_2,\lambda_{3},\lambda_{112})$ be the quotient of $T(V)$ modulo the ideal generated by
	\begin{align*}
		y_{1}^3 &= \lambda_1; &  y_{2}^{6} &=\lambda_2; &  y_{221} &=\lambda_{3};
	\end{align*}
	\vsp
	\begin{align*}
		y_{112}^{2}+4q_{21}\zeta\lambda_1y_2y_{12} &=\lambda_{112}.
	\end{align*}
	Then $\mE\in\hmod$ and $\mA(\lambda_1,\lambda_2,\lambda_{3},\lambda_{112})\coloneqq \mE\#H \in\Cleft\B(V)\#H$.
\end{prop}

\pf
See \texttt{cleft-br2a-q-z.(g|log)}.
Here $\lambda_1\neq 0$ only if $q_{21}^3=1$ and $\lambda_2\neq 0$ only if $q_{12}^6=1$. As well, $\lambda_{3}\neq 0$ only if 
$\zeta=q_{21}=-q_{12}$. In this case, the previous conditions hold and hence we have that when $\bq=\left(\begin{smallmatrix}\zeta& -\zeta \\ \zeta & -\zeta
\end{smallmatrix}\right)$, we can have all deformation parameters $\lambda_1, \lambda_2,\lambda_{3}$ and $\lambda_{112}$ nonzero, simultaneously.
For any other braiding matrix $\lambda_{3}=0$.
\epf

\begin{prop}\label{pro:lifting-br(2,a)-a-q-zeta}
	Fix $V$ of type
	$\br(2, a)$,
	as in \emph{(\ref{eq:dynkin-br(2,a)} a)}, with $q=-\zeta$. Let $\lambda_1, \lambda_2,\lambda_{3},\lambda_{112}\in\k$ subject to \eqref{eqn:mu} and let $\mL=\mL(\lambda_1,\lambda_2,\lambda_{3},\lambda_{112})$ be the quotient $T(V)\# H / \ideal$, where
	$\ideal$ is the ideal of $T(V)\# H$ generated by
	\begin{align*}
		a_{1}^3 &= \lambda_1(1-g_1^3); &  a_{2}^{6} &=\lambda_2(1-g_2^6); &  a_{221} &=\lambda_{3}(1-g_{122});
	\end{align*}
	\vsp
	\begin{align*}
		a_{112}^{2} &=\lambda_{112}(1-g_{112}^2)-s_{112},
	\end{align*}
	\noindent 
	where
	\[
	s_{112}=4\zeta q_{21}\lambda_1a_2a_{12} -4\zeta\lambda_1\lambda_{3}g_{122}(1-g_1^3).
	\] 
	Then $\mL$ is a Hopf algebra and a cocycle deformation of $\B(V)\# H$ such that $\gr\mL\simeq \B(V)\# H$.
	
	Conversely, if $L$ is a Hopf algebra such that its infinitesimal braiding is a principal realization $V\in\ydh$, then there are $\lambda_1, \lambda_2,\lambda_{3},\lambda_{112}\in\k$ such that 
	$L\simeq \mL(\lambda_1, \lambda_2,\lambda_{3},\lambda_{112})$.
\end{prop}

\pf 
See \texttt{lift-br2a-q-z.(g|log)}
\epf

\subsection{The case $q\neq -1,-\zeta$}\label{subsubsec:br(2,a)-q-generico}
Let us fix 
\begin{align*}
	M := \ord \zeta q^{-1}.
\end{align*} 
The Nichols algebra $\toba(V)$ is
generated by $x_1, x_{2}$
with defining relations
\begin{align*}
	x_{1}^3 &= 0; & x_{2}^{N} &=0; & x_{221} &=0;
\end{align*}
\vsp
\begin{align*}
	x_{112}^{M} &=0.
\end{align*}

\begin{rem}
	We see that, depending on $N$, we may have
	\begin{align*}
		M&=N, &\text{ or }&& M&=3N, &\text{ or }&& M&=N/3.
	\end{align*}
	In any case, we observe that relation $x_{221}=0$ is not deformed, as $\chi_1\chi_2^2\neq\eps$.
\end{rem}

We will compute the cleft objects and liftings for a pair $(N,M)$ corresponding to each case, namely we will work with
\begin{enumerate}
	\item $N=4$, so $M=3N=12$.
	\item $N=6$, so $M=N=6$ (here we necessarily have $q=-\zeta^2$).
	\item $N=12$, so $M=N/3=4$.
\end{enumerate}

\subsubsection{Case $(N,M)=(4,12)$}

\begin{prop}\label{pro:cleft-br(2,a)-q-generico-N4}
	Fix $V$ of type
	$\br(2, a)$,
	as in \emph{(\ref{eq:dynkin-br(2,a)} a)}, with $q\neq-1,-\zeta$. Let $\lambda_1,\lambda_2,\lambda_{112}\in\k$ subject to \eqref{eqn:mu} and let $\mE=\mE(\lambda_1,\lambda_2,\lambda_{112})$ be the quotient of $T(V)$ modulo the ideal generated by
	\begin{align*}
		y_{1}^3 &=\lambda_1; &  y_{2}^{N} &=\lambda_2; &  y_{221} &=0; &
		y_{112}^{M} &=\lambda_{112}.
	\end{align*}
	Then $\mE\in\hmod$ and $\mA(\lambda_1,\lambda_2,\lambda_{112})\coloneqq \mE\#H \in\Cleft\B(V)\#H$.
\end{prop}
\pf
See \texttt{cleft-br2a-M12-N4.(g|log)}
\epf

\begin{prop}\label{pro:lifting-br(2,a)-a-q-generico-N4}
	Fix $V$ of type
	$\br(2, a)$,
	as in \emph{(\ref{eq:dynkin-br(2,a)} a)}, with $q\neq-1,-\zeta$. Let $\lambda_1,\lambda_2,\lambda_{112}\in\k$ subject to \eqref{eqn:mu} and let $\mL=\mL(\lambda_1,\lambda_2,\lambda_{112})$ be the quotient $T(V)\# H / \ideal$, where
	$\ideal$ is the ideal of $T(V)\# H$ generated by
	\begin{align}
		a_{1}^3 &= \lambda_1(1-g_1^3); & a_{2}^{N} &=\lambda_2(1-g_2^N); &  a_{221} &=0;
	\end{align}
	\vsp
	\begin{align*}
		a_{112}^{M} &=\lambda_{112}(1-g_{112}^M).
	\end{align*}
	Then $\mL$ is a Hopf algebra and a cocycle deformation of $\B(V)\# H$ such that $\gr\mL\simeq \B(V)\# H$.
	
	Conversely, if $L$ is a Hopf algebra such that its infinitesimal braiding is a principal realization $V\in\ydh$, then there are $\lambda_1, \lambda_2,\lambda_{112}\in\k$ such that 
	$L\simeq \mL(\lambda_1, \lambda_2,\lambda_{112})$.\end{prop}
\pf
See \verb!lift-br2a-M12-N4-with_c2.(g|log)!
for the case $\lambda_2\neq 0$ (so $\lambda_1=0$) and 
\verb!lift-br2a-M12-N4-t1-PBW.(g|log)!
for $\lambda_1\neq 0$, $\lambda_2=0$.
\epf
\begin{rem}
	In the \texttt{GAP} files used in proof of Proposition \ref{pro:lifting-br(2,a)-a-q-generico-N4}, we have used new ideas to compute the liftings. When $\lambda_2\neq 0$, we introduce a new variable $c_2=\lambda_2(1-g_2^4)$ and use that it commutes with $x_1$ and $x_2$. On the other hand, this is no longer true when $\lambda_1\neq 0$, since $\lambda_1(1-g_1^3)$ does not commute with $x_1$. For this case, we have used the whole PBW basis in the letters $x_2,x_{112},x_{12},x_1$ inherited from the Nichols algebra $\B(V)$ instead of the the letters $x_2,x_1$ usually used to deal with the liftings.
\end{rem}

\subsubsection{Case $(N,M)=(6,6)$}

\begin{prop}\label{pro:cleft-br(2,a)-q-generico-N6}
	Fix $V$ of type
	$\br(2, a)$,
	as in \emph{(\ref{eq:dynkin-br(2,a)} a)}, with $q\neq-1,-\zeta$. Let $\lambda_1,\lambda_2,\lambda_{112}\in\k$ subject to \eqref{eqn:mu} and let $\mE=\mE(\lambda_1,\lambda_2,\lambda_{112})$ be the quotient of $T(V)$ modulo the ideal generated by
	\begin{align*}
		y_{1}^3 &=\lambda_1; &  y_{2}^{6} &=\lambda_2; &  y_{221} &=0; &
		y_{112}^{6} &=\lambda_{112}.
	\end{align*}
	Then $\mE\in\hmod$ and $\mA(\lambda_1,\lambda_2,\lambda_{112})\coloneqq \mE\#H \in\Cleft\B(V)\#H$.
\end{prop}

\pf
See \texttt{cleft-br2a-M-N-6.(g|log)}.
\epf

\begin{prop}\label{pro:lifting-br(2,a)-a-q-generico-N6}
	Fix $V$ of type $\br(2, a)$,
	as in \emph{(\ref{eq:dynkin-br(2,a)} a)}, with $q\neq-1,-\zeta$. Let $\lambda_1,\lambda_2,\lambda_{112}\in\k$ subject to \eqref{eqn:mu} and let $\mL=\mL(\lambda_1,\lambda_2,\lambda_{112})$ be the quotient $T(V)\# H / \ideal$, where
	$\ideal$ is the ideal of $T(V)\# H$ generated by
	\begin{align*}
		a_{1}^3 &= \lambda_1(1-g_1^3); & a_{2}^{6} &=\lambda_2(1-g_2^6); &  a_{221} &=0;
	\end{align*}
	\vsp
	\begin{align*}
		a_{112}^{6} &=\lambda_{112}(1-g_{112}^6)-s_{112},
	\end{align*}
	\noindent 
	where
	\begin{align*}
		s_{112}= \lambda_1^4\lambda_2g_2^6(1-g_1^{12}).
	\end{align*}
	Then $\mL$ is a Hopf algebra and a cocycle deformation of $\B(V)\# H$ such that $\gr\mL\simeq \B(V)\# H$.
	
	Conversely, if $L$ is a Hopf algebra such that its infinitesimal braiding is a principal realization $V\in\ydh$, then there are $\lambda_1, \lambda_2,\lambda_{112}\in\k$ such that 
	$L\simeq \mL(\lambda_1, \lambda_2,\lambda_{112})$.
\end{prop}
\pf
See \texttt{lift-br2a-M-N-6.(g|log)}
\epf

\subsubsection{Case $(N,M)=(12,4)$}

\begin{prop}\label{pro:cleft-br(2,a)-q-generico-N12}
	Fix $V$ of type
	$\br(2, a)$,
	as in \emph{(\ref{eq:dynkin-br(2,a)} a)}, with $q\neq-1,-\zeta$. Let $\lambda_1,\lambda_2,\lambda_{112}\in\k$ subject to \eqref{eqn:mu} and let $\mE=\mE(\lambda_1,\lambda_2,\lambda_{112})$ be the quotient of $T(V)$ modulo the ideal generated by
	\begin{align*}
		y_{1}^3 &=\lambda_1; &  y_{2}^{N} &=\lambda_2; &  y_{221} &=0; &
		y_{112}^{M}-s_{112}=\lambda_{112}.
	\end{align*}
	\noindent where
	\[
	s_{112}\coloneqq (4q^4+8q^7-4q^8+8q^{11})q^2q_{21}\lambda_1^2y_2^2y_{12}^2.
	\]
	Then $\mE\in\hmod$ and $\mA(\lambda_1,\lambda_2,\lambda_{112})\coloneqq \mE\#H \in\Cleft\B(V)\#H$.
\end{prop}

\pf
See \texttt{cleft-br2a-M4-N12-E(12).(g|log)}.
\epf

\begin{prop}\label{pro:lifting-br(2,a)-a-q-generico-N12}
	Fix $V$ of type
	$\br(2, a)$,
	as in \emph{(\ref{eq:dynkin-br(2,a)} a)}, with $q\neq-1,-\zeta$.Let $\lambda_1, \lambda_2,\lambda_{112}\in\k$ subject to \eqref{eqn:mu} and let $\mL=\mL(\lambda_1,\lambda_2,\lambda_{112})$ be the quotient $T(V)\# H / \ideal$, where
	$\ideal$ is the ideal of $T(V)\# H$ generated by
	\begin{align}\label{eq:rels-deformadas-br(2,a)-q-generico}
		a_{1}^3 &= \lambda_1(1-g_1^3); & a_{2}^{N} &=\lambda_2(1-g_2^N); &  a_{221} &=0;
	\end{align}
	\vsp
	\begin{align*}
		a_{112}^{M} &=\lambda_{112}(1-g_{112}^M)-s_{112},
	\end{align*}
	\noindent where
	\begin{align*}
		s_{112}\coloneqq & 
		-4q^4(2+q^4+2q^{7})\lambda_1^2q_{21}a_2^2a_{12}^2.
	\end{align*}
	Then $\mL$ is a Hopf algebra and a cocycle deformation of $\B(V)\# H$ such that $\gr\mL\simeq \B(V)\# H$.
	
	Conversely, if $L$ is a Hopf algebra such that its infinitesimal braiding is a principal realization $V\in\ydh$, then there are $\lambda_1, \lambda_2,\lambda_{112}\in\k$ such that 
	$L\simeq \mL(\lambda_1, \lambda_2,\lambda_{112})$.
\end{prop}

\pf
See \texttt{lift-br2a-M4-N12-E(12)-PBW.(g|log)}.
\epf

\end{document}